\documentclass{article}
\usepackage{amsmath,amssymb,amsthm}
\usepackage{graphicx}
\usepackage[all]{xy}
\theoremstyle{plain}

\title{Erratum: Linear projections and successive minima}
\author{Christophe {\sc Soul\'e}\footnote{CNRS et IH\'ES, Le Bois-Marie, 35 route de Chartres, 91440 {\sc Bures-sur-Yvette}, France \newline soule@ihes.fr}}

\begin{document}

\maketitle  

\section{Erratum}

The proof of Proposition~1 and Theorem~2 in  \cite{So} is incorrect. Indeed, \S2.5 and \S2.7 in {\it op.cit} contain a vicious circle: the definition of the filtration $V_i$, $1 \leq i \leq n$, in \S2.5 depends on the choice of the integers $n_i$, when the definition of the integers $n_i$ in \S2.7 depends on the choice of the filtration $(V_i)$. Thus, only Theorem~1 and Corollary~1 in \cite{So} are proved.
We shall prove below another result instead of Proposition~1 in \cite{So}.

\smallskip

I thank J.-B. Bost, C. Gasbarri and C. Voisin for their help.

\section{An inequality}

\subsection{ \ }

Let $K$ be a number field, $O_K$ its ring of algebraic integers and $S = {\rm Spec} (O_K)$ the associated scheme. Consider an hermitian vector bundle $(E,h)$ over $S$. Define the $i$-th successive minima $\mu_i$ of $(E,h)$ as in \cite{So} \S2.1. Let $X_K \subset {\mathbb P} (E_K^{\vee})$ be a smooth, geometrically irreducible curve of genus $g$ and degree $d$. We assume that $X_K \subset {\mathbb P} (E_K^{\vee})$ is defined by a complete linear series on $X_K$ and that $d \geq 2g+1$. The rank of $E$ is thus $N = d+1-g$. Let $h (X_K)$ be the Faltings height of $X_K$ (\cite{So} \S2.2).

\smallskip

For any positive integer $i \leq N$ we define the integer $f_i$ by the formulas
$$
\begin{matrix}
f_i = i-1 \hfill &\mbox{if} &i-1 \leq d - 2g \, , \hfill \\
f_i = i-1+\alpha &\mbox{if} &i-1 = d-2g+\alpha \, , \quad 0 \leq \alpha \leq g \, ,
\end{matrix}
$$
and $f_N = d$.

\smallskip

Fix two natural integers $s$ and $t$ and suppose that $1 \leq s < t \leq N-2$. When $2 \leq i \leq s$ we let
$$
A_i  = \frac{f_i^2}{(i-1) \, f_i - \underset{j=2}{\overset{i-1}{\sum}} \, f_j} \, ,
$$
and, when $t \leq i \leq N$,
$$A_i  = \frac{f_i^2}{((i-t+s) \, f_i - (f_1+f_2+ \ldots +f_s+f_t+ \ldots +f_{i-1}))} \, ,
$$
consider
$$A(s,t) = \max_{2 \leq i \leq s \, {\rm or} \, t\leq i \leq N} A_i \, .
$$

\bigskip

\noindent {\bf Theorem 1.} {\it  There exists a constant $c(d)$ such that the following inequality holds:}
$$
\frac{h(X_K)}{[K:{\mathbb Q}]} + (2d - A(s,t) (N-t+s+1)) \, \mu_1 + A(s,t) (\sum_{\alpha = 1}^{N+1-t}\mu_{\alpha} + \sum_{\alpha = N+1-s}^{N}\mu_{\alpha})  + c(d) \geq 0 \, .
$$

\subsection{ \ }

To prove Theorem~1 we start by the following variant of Corollary 1 in \cite{MO}.

\bigskip

\noindent {\bf Proposition 1.} {\it  Fix an increasing sequence of integers $0 = e_1 \leq e_2 \leq \ldots \leq e_N$ and a decreasing sequence of numbers $r_1 \geq r_2 \geq \ldots \geq r_N$. Assume that $e_s = e_{s+1} = \ldots = e_{t-1} $. Let
$$
S = \min_{0 = i_0 < \ldots < i_{\ell} = N} \, \sum_{j=0}^{\ell - 1} \, (r_{i_j} - r_{i_{j+1}})(e_{i_j} + e_{i_{j+1}}) \, .
$$
Then}
$$
S \leq B(s,t) ( \sum_{j=1}^{s} (r_j - r_N) + \sum_{j=t}^{N} (r_j - r_N))  \, ,
$$
{\it where
$$
B(s,t) = \max_{2 \leq i \leq s \, {\rm or} \, t\leq i \leq N} B_i \, ,
$$
and $B_i$ is defined by the same formula as $A_i$, each $f_j$ being replaced by $e_j$. 
\bigskip}

\noindent {\bf Proof.} We can assume that $r_N=0$. As in \cite{MO}, proof of Theorem 1, we may first assume that $S=1$ and seek to minimize $\underset{j=1}{\overset{s}{\sum}} \ r_j + \underset{j=t}{\overset{N}{\sum}} \ r_j $. If we graph the points $(e_j , r_j)$, $S/2$ is the area of the Newton polygon they determine in the first quadrant. Moving the points not lying on the polygon down onto it only reduces $\underset{j=1}{\overset{s}{\sum}} \ r_j + \underset{j=t}{\overset{N}{\sum}} \ r_j$, so we may assume that all the points actually lie on the polygon. In particular 
$$
S(r_1 , \ldots , r_N) =
S(r_1 , r_2,\ldots ,r_{s} , \ldots , r_s, r_{t}, \ldots, r_N) $$ and we may assume that the point $(e_j , r_j) = (e_s,r_j)$ lies on this polygon when $s \leq j \leq t-1$. For such $r_i$'s we have
$$
S(r_1 , \ldots , r_N) = \sum_{i=1}^{N-1} (r_i - r_{i+1}) \, (e_i + e_{i+1}) \, .
$$

Let $\sigma_i = r_{i-1} - r_i$, $i= 2,\ldots , N$. The condition that the points $(e_i , r_i)$ lie on their Newton polygon and that the $r_i$ decrease becomes, in terms of the $\sigma_i$,
\begin{equation}
\label{1}
\frac{\sigma_{2}}{e_{2} - e_{1}} \geq \frac{\sigma_{3}}{e_{3} - e_{2}} \geq \ldots \geq 0 \, .
\end{equation}
Furthermore
$$
 \sigma_{s+1} = \ldots = \sigma_{t-1} = 0 \, .
$$
Next, we impose the constraint $\underset{j=1}{\overset{s}{\sum}} \ r_j + \underset{j=t}{\overset{N}{\sum}} \ r_j = 1$, i.e.
\begin{equation}
\label{2}
\sum_{j=2}^s  \, (j-1) \, \sigma_j + \sum_{j=t}^N (j-t+s) \, \sigma_j = 1 \, .
\end{equation}
In the  subspace of the points $\sigma = (\sigma_{2} , \ldots , \sigma_s, \sigma_t,\ldots,\sigma_N)$ defined by (\ref{2}), the inequalities (\ref{1}) define a simplex. The linear function
$$
S = \sum_{2 \leq j \leq s} \sigma_j \, (e_{j-1} + e_j)+\sum_{ t \leq j \leq N} \sigma_j \, (e_{j-1} + e_j)
$$
must achieve its maximum on this simplex at one of the vertices, {\it i.e.} a point where, for some $i$ and $\alpha$, we have
$$
\alpha = \frac{\sigma_{2}}{e_{2} - e_{1}} = \ldots = \frac{\sigma_i}{e_i - e_{i-1}} > \frac{\sigma_{i+1}}{e_{i+1} - e_i} = \ldots = 0 \, .
$$
We get
$$
\sigma_j = \left\{
\begin{matrix}
\alpha (e_j - e_{j-1}) &\mbox{if} \, j \leq i \\
0 \hfill &\mbox{else.} \hfill
\end{matrix} 
\right.
$$
Then, using (\ref{2}) we get, if $i \leq s$,
$$
\alpha = \left( (i-1)\, e_i - \sum^{i-1}_{j = 2} e_j \right)^{-1} \, ,
$$
and, when $i\geq t$,
$$
\alpha = \left( (i-t+s) \, e_i - e_1 - e_2 - \ldots -e_s - e_t - \ldots - e_{i-1}\right)^{-1} \, .
$$
Since
$$ S = \alpha \, \sum^i_{j=2} (e_j^2 - e_{j-1}^2)  = \alpha \, e_i^2 $$
 Proposition 1  follows.
\bigskip

\subsection{ \ }

We come back to the situation of Theorem~1. For every complex embedding $\sigma : K \to {\mathbb C}$, the metric $h$ defines a scalar product $h_{\sigma}$ on $E \underset{O_K}{\otimes} {\mathbb C}$. If $v \in E$ we let
$$
\Vert v \Vert = \max_{\sigma} \, \sqrt{h_{\sigma} (v,v)} \, .
$$
Choose $N$ elements $x_1 , \ldots , x_N$ in $E$, linearly independent over $K$ and such that
$$
\log \Vert x_i \Vert = \mu_{N-i+1} \, , \quad 1 \leq i \leq N \, .
$$
Let $y_1 , \ldots , y_N \in E_K^{\vee}$ be the dual basis of $x_1 , \ldots , x_N$. Let $A(d)$ be the constant appearing in \cite{So}, Theorem~1. From \cite{So}, Corollary~1, we deduce

\bigskip

\noindent {\bf Lemma 1.} {\it Assume $1 \leq s \leq t \leq N-2$. We may choose integers $n_i$, $s+1\leq i \leq t-1$, such that}
\begin{enumerate}
\item[i)] {\it For all $i$ \, $\vert n_i \vert \leq A(d) + d$}
\item[ii)] {\it Let $w_i = y_i$ if $1 \leq i \leq s$ or $t\leq i \leq N$
and
$w_i = y_i + n_i \, y_{i+1} \, {\it if} \, s+1 \leq i \leq t-1 .$
Let $\langle w_1 , \ldots , w_i \rangle \subset E_K^{\vee}$ be the subspace spanned by $w_1 , \ldots , w_i$, and
$$
W_i = E_K^{\vee} / \langle w_1 , \ldots , w_i \rangle \, 
$$
($W_0 = E_K^{\vee}$).
Then, when $s+1 \leq i \leq t-1$, the linear projection from ${\mathbb P} (W_{i-1})$ to ${\mathbb P} (W_i)$ does not change the degree of the image of $X_K$.}
\end{enumerate}
\bigskip

\subsection{ \ }

Let $(v_i) \in E_K^N$ be the dual basis of $(w_i)$. We have
$$ v_i = x_i \, {\rm when} \, i\leq s+1 \, {\rm or} \, i\geq t+1$$
and
\begin{equation}
v_i = x_i -n_{i-1} \, x_{i-1} + n_{i-1} \, n_{i-2} \, x_{i-2} - \ldots \pm n_{i-1} \ldots n_{s+1} \, x_{s+1} \, \nonumber
\end{equation}
when $s+2 \leq i \leq t$.

From these formulas it follows that 
 there exists a positive constant $c_1(d)$ such that

$$
\log \Vert v_i \Vert \leq r_i = \left\{
\begin{matrix}
\mu_{N+1-i} + c_1 (d) &\mbox{if} \, i\leq s  \, \mbox{or} \, i \geq t+1  \hfill \\
\mu_{N-s} + c_1 (d) &\mbox{if} \, s+1 \leq i \leq t \, .
\end{matrix}
\right.
$$
Let $d_i$ be the degree  of the image of $X_K$ in ${\mathbb P} (W_i)$, and $e_i = d - d_i$.
By Lemma~1 we have
$$
e_s = e_{s+1} = \ldots = e_{t-1}\, .$$
Therefore we can argue
 as in \cite{S1}, Theorem 1 and \cite{So} pp.~50--53, to deduce Theorem~1 from Proposition~1.
\vglue 2cm


\begin{thebibliography}{99}
\bibitem{MO}Morrison, I.
Projective stability of ruled surfaces. 
Invent. Math. 56, 269-304 (1980).
\bibitem{S1}Soul\'e, C. Successive minima on arithmetic varieties. Compositio Mathematica 96 (1995), 85-98.
\bibitem{So}Soul\'e, C. Linear projections and successive minima. Nagoya Math. J. 197 (2010), 45-57.
\end{thebibliography}
\end{document}